\RequirePackage{fix-cm}
\documentclass{svjour3-1}                     % onecolumn (standard format)
\smartqed  % flush right qed marks, e.g. at end of proof

\usepackage{amssymb,amsmath,url}
\usepackage{graphicx,latexsym,amssymb,color}

\newtheorem{thm}{Theorem}

\newtheorem{lem}{Lemma}
\newtheorem{cor}{Corollary}

\newtheorem{defn}{Definition}

\begin{document}

\title{A general lower bound for the domination number of cylindrical graphs
%\thanks{Partially supported by grants TIN2015-66680 (MINECO-ERDF),  MTM2015-63791-R (MINECO/FEDER) and FQM305 (Junta de Andaluc\'ia). }
}


\author{J.J. Carre\~{n}o,  J.A. Mart\'inez, M.L. Puertas}
%Jos\'e Juan Carre\~{n}o,  Jos\'e Antonio Mart\'inez, Mar\'ia Luz Puertas

\authorrunning{J.J. Carre\~{n}o \and J. A. Mart\'inez \and M.L. Puertas} % if too long for running head

\institute{J.J. Carre\~{n}o\at
              Department of Applied Mathematics for Information and Communication Technologies\\
Universidad Polit\'ecnica de Madrid, Spain\\
           J.A. Mart\'inez\at
              Department of Computer Science, Universidad de Almer\'ia, Spain\\
          M.L. Puertas \at
              Department of Mathematics, Universidad de Almer\'ia, Spain
              \and
              %\email{jjcc@etsisi.upm.es, jmartine@ual.es, mpuertas@ual.es}
}


\maketitle

\vspace{-2cm}
\begin{abstract}
In this paper we present a lower bound for the domination number of the Cartesian product of a path and a cycle, that is tight if the length of the cycle is a multiple of five. This bound improves the natural lower bound obtained by using the domination number of the Cartesian product of two paths, that is the best one known so far.
%\subclass{05C38\and 05C76\and 05C85}
\end{abstract}

\section{Introduction}\label{sec:introduction}

A \emph{dominating set} in a graph $G$ is a vertex subset $S$ such that every vertex in $V(G)\setminus S$
has a neighbor in $S$. The \emph{domination number} of $G$ is the minimum size of a dominating set of $G$ and it is denoted by $\gamma(G)$.

The domination number is a classical graph parameter that has a challenging de\-ve\-lop\-ment in Cartesian product graphs. As proof of this, the computation of the domination number of the grid graphs, that is, the Cartesian product of two paths, was an open problem for almost thirty years. The final paper concerning this problem~\cite{Goncalves2011} closes a list of partial results regarding exact values for particular cases~\cite{Alanko2011,Chang1993,Cockayne1985,Jacobson1983} and general upper and lower bounds~\cite{Chang1992,Cockayne1985,Guichar2003}, among other papers.

Contrary to grids, the complete computation of the domination number of cylindrical graphs, that is, the Cartesian product of a path and a cycle, is still open. Currently, there are some partial results, see for instance~\cite{Crevals2014,Klavzar1996,Pavlic2013}, showing different techniques to compute exact values of particular cases, where one of the parameters, the path length or the cycle length, is small enough.

We devote this paper to computing a general lower bound for the domination number of the cylindrical graphs, that is tight if the length of the cycle is a multiple of five. Our strategy to obtain this bound is similar to the approach in~\cite{Guichar2003}, that uses the concept of wasted domination. On the other hand, we found inspiration in~\cite{Klavzar1996} to design a $(\min,+)$ matrix multiplication algorithm that computes the minimum wasted domination needed for our purpose.

For undefined concepts about graphs and the Cartesian product see~\cite{Chartrand2011,Imrich2000}. Let $G$ be a graph with vertex set $V(G)$ and edge set $E(G)$, the \emph{closed neighbor} of a vertex $u\in V(G)$ is $N[u]=\{u\}\cup \{ v\in V(G)\colon uv\in E(G)\}$. The closed neighbor of a vertex subset $U\subseteq V(G)$ is $N[U]=\bigcup \limits_{u\in U} N[u]$.

\begin{figure}[hbp]
\begin{center}
\includegraphics[width=.4\textwidth]{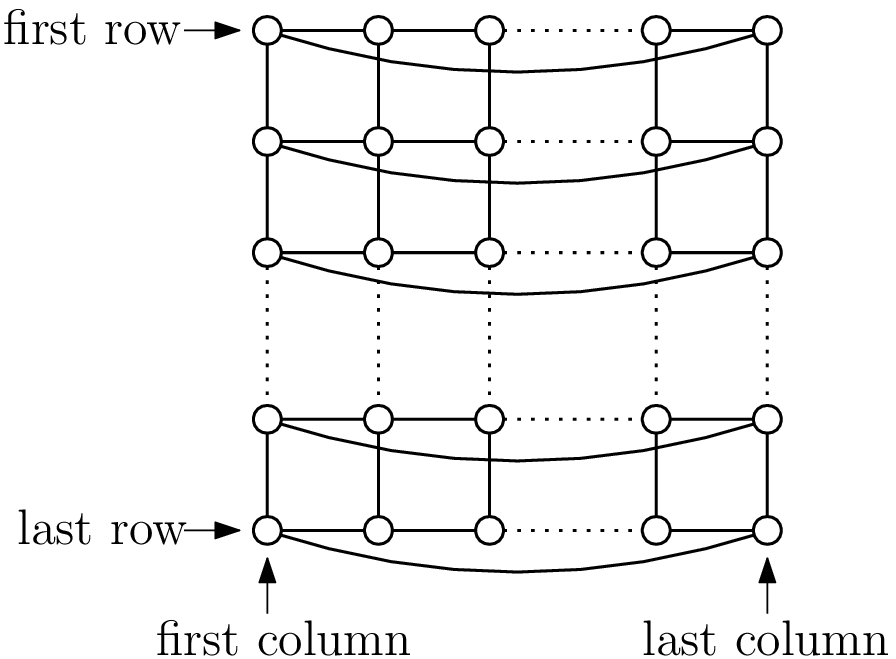}
\caption{The cylinder $P_m \Box C_n$}
\label{fig:cylinder}
\end{center}
\end{figure}

For integers $m\geq 2$ and $n\geq 3$, the Cartesian product of the path $P_m$ and the cycle $C_n$ is the cylinder $P_m\Box C_n$, that has vertex set
$$V(P_m\Box C_n)=\{v_{ij}\colon i\in \{0, \dots, m-1 \}, j\in \{0,\dots , n-1\} \}$$
and there is an edge between two different vertices $v_{ij}$ and $v_{i'j'}$ if and only if they satisfy one of the following properties:\\
$\begin{array}{l}
\bullet \ i=i', j=0, j'=n-1,\\
\bullet \ i=i', j=n-1, j'=0,\\
\bullet \ i=i', |i-i'|=1,\\
\bullet \ |i-i'|=1, j=j'.
\end{array}
$\\
We say that $P_m\Box C_n$ has $m$ rows and $n$ columns, each row is a cycle with $n$ vertices and each column is a path with $m$ vertices (see Figure~\ref{fig:cylinder}). We numerate rows from top to bottom and columns from left to right. Moreover, we consider that the last column is the previous column of the first one (and the first column follows the last one).

The paper is organized as follows: in Section~\ref{sec:known} we summarize the known values and bounds for the domination number of cylindrical graphs. In Section~\ref{sec:wasted} we recall the concept of wasted domination, which goes on to play a central role in the computation of our lower bound, that is presented in Theorem~\ref{thm:bound}, in Section~\ref{sec:bound}. Conclusions are finally shown in Section~\ref{sec:conclusions}.

\section{Known bounds and exact values}\label{sec:known}

Within this section we summarize the exact values computed up to now and the known lower and upper bounds for $\gamma(P_m\Box C_n)$. In~\cite{Crevals2014}, the author provides formulas for $\gamma(P_m\Box C_n)$, in cases $2\leq m\leq 22$ and $n\geq m$ and also in cases $3\leq n\leq 30$ and $m\geq n$. Several other papers exist showing particular small cases~\cite{Nandi2011,Pavlic2013}, all of which are included in~\cite{Crevals2014}.

Since the domination number of $P_m\Box C_n$ is known when $m$ or $n$ is small, lower and upper bounds are interesting when both are big enough. Regarding the known bounds, some can be deduced from the relationship between cylinders and grids, or cylinders and tori, while others come from specific constructions. We begin with the second type. Particular dominating sets of $P_m\Box C_n$, with $m\geq n$, based on regular patterns, give that (see~\cite{Pavlic2013})
\begin{equation}\label{eq:construction}
\gamma(P_m\Box C_n)\leq
\left\{
\def\arraystretch{2}
\displaystyle
\begin{array}{ll}
(m+2)k=(m+2)\displaystyle\frac{n}{5}& \text{if } n=5k,\\
(m+2)\displaystyle\frac{8k+3}{8}=(m+2)\displaystyle\Big(\frac{n}{5}+\frac{7}{40}\Big) & \text{if } n=5k+1,\\
(m+2)\displaystyle\frac{2k+1}{2}=(m+2)\displaystyle\Big(\frac{n}{5}+\frac{1}{10}\Big)  & \text{if } n=5k+2,\\
(m+2)(k+1)=(m+2)\displaystyle\Big(\frac{n}{5}+\frac{2}{5}\Big)  & \text{if } n=5k+3,\\
(m+2)(k+1)=(m+2)\displaystyle\Big(\frac{n}{5}+\frac{1}{5}\Big)  & \text{if } n=5k+4.\\
\end{array}
\right.
\end{equation}
Authors conjecture that the bounds are exact values for $n = 5k, 5k+1, 5k+2, 5k+4$. Unfortunately, values computed in~\cite{Crevals2014} disprove this conjecture for $n=5k+1,5k+2$, $5k+4$, (and moreover for $n=5k+3$). For instance, for $m=22$:

$
\def\arraystretch{1.5}
\begin{array}{ll}
n=18, (k=3): & \gamma(P_{22}\Box C_{18})=89<(22+2)(3+1)=96,\\
n=19, (k=3): & \gamma(P_{22}\Box C_{19})=93<(22+2)(3+1)=96,\\
n=21, (k=4): & \gamma(P_{22}\Box C_{21})=102< (22+2)\frac{(8\cdot 4 +3)}{8}=105,\\
n=22, (k=4): & \gamma(P_{22}\Box C_{22})=107< (22+2)\frac{(2\cdot 4 +1)}{2}=108.\\
\end{array}
$
\par\bigskip

\begin{figure}[htp]
\begin{center}
\includegraphics[width=.45\textwidth]{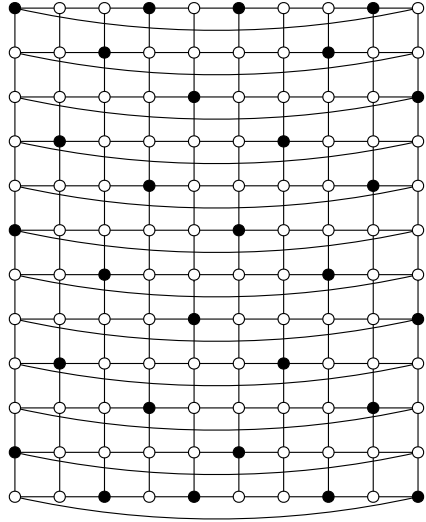}
\caption{Dominating set of $P_{12} \Box C_{10}$ with $\frac{(12+2)10}{5}$ vertices and following a regular pattern}
\label{fig:5-cycle}
\end{center}
\end{figure}

However, the comparison with the exact values given in~\cite{Crevals2014}, for the cases $16\leq m\leq 22$, leads us to think that the construction shown in~\cite{Pavlic2013} for the case $n=5k$ is optimal, for $m$ big enough (see Figure~\ref{fig:5-cycle}). Clearly this construction can be done for $5\leq n\equiv 0\pmod 5$ and any $m\geq 2$ and it provides a dominating set with $\frac{(m+2)n}{5}=(m+2)k$ vertices.

With reference to bounds obtained from the relationship of cylinders with grids and tori, we will consider that $P_m\Box P_n$, $P_m\Box C_n$ and $C_m\Box C_n$ have the same vertex set. The domination number of grids was computed in \cite{Goncalves2011}, where authors obtained
$$\gamma(P_m\Box P_n)=\Big\lfloor \displaystyle\frac{(m+2)(n+2)}{5}\Big\rfloor -4, \text{ for }  n,m\geq 16.$$

It is clear that a dominating set of $P_m\Box P_n$ also dominates the cylinder of the same size $P_m\Box C_n$, therefore
\begin{equation}\label{eq:upper}
\gamma(P_m\Box C_n)\leq \Big\lfloor \displaystyle\frac{(m+2)(n+2)}{5}\Big\rfloor -4, \text{ for }  n,m\geq 16.
\end{equation}

Note that this bound does not improve, in general, those given in Equation~\ref{eq:construction} (just in the case $n=5k+3$ the bound in Equation~\ref{eq:upper} is better for every $m\geq 16$).

On the other hand, let $S\subseteq V(P_m\Box C_n)$ be a minimum dominating set, and let $P_m\Box P_{n-2}$ be the subgraph of $P_m\Box P_n$ avoiding the first and the last columns. We construct the set $S'$ by adding to $S \cap V(P_m\Box P_{n-2})$, the neighbors in $V(P_m\Box P_{n-2})$ of the vertices in $S$ in the first and the last column of $P_m\Box P_n$. Clearly $S'$ is a dominating set of $P_m\Box P_{n-2}$ and it satisfies $|S'|\leq |S|$, therefore, for $n,m\geq 18$:
\begin{equation}\label{eq:lower}
\Big\lfloor \displaystyle\frac{(m+2)n}{5}\Big\rfloor -4=\gamma(P_m\Box P_{n-2})\leq |S'|\leq |S|=\gamma(P_m\Box C_n).
\end{equation}
A dominating set of $P_m\Box C_n$ also dominates $C_m\Box C_n$, so $\gamma(C_m\Box C_n)\leq \gamma(P_m\Box C_n)$. The domination number of torus $C_m\Box C_n$ remains unknown in the general case. Exact values for $3\leq m\leq 20$ and $m\leq n$ are obtained in \cite{Crevals2018}, while general lower bounds can be found in \cite{El-Zahar2002,VanWieren2007}. These lower bounds, that depend of the parity of $n$ module $5$, are smaller than the lower bound given in Equation~\ref{eq:lower}, so they do not provide an improved bound.

To sum up, we can locate the domination number of the cylinder between the following two bounds, for big enough $m\geq n=5k+r$:
\begin{equation}\label{eg:bothbounds}
\displaystyle\Big\lfloor \frac{(m+2)n}{5}\Big\rfloor  -4\leq \gamma(P_m\Box C_n)\leq
\left\{
\def\arraystretch{1.5}
\begin{array}{ll}
(m+2)k & \text{ if } n\equiv 0\!\!\!\!\!\pmod 5,\\
(m+2)(k+1) & \text{ otherwise. }\\
\end{array}
\right.
\end{equation}

In the cases $n\equiv 1,2\pmod 5$, we have rounded up the upper bound given in Equation~\ref{eq:construction}, that we know is not tight, in order to obtain a regular formula.

In the case $n\equiv 0\pmod 5$, the construction provided in~\cite{Pavlic2013}, that we have shown in Figure~\ref{fig:5-cycle}, of a dominating set of $P_m\Box C_{5k}$ with $(m+2)k$ vertices, is likely to be the optimal construction. If this is true, the lower bound of Equation~\ref{eg:bothbounds} is not tight in this case and following this intuition, we have conjectured how a better general lower bound should be. It has led us to prove in Theorem~\ref{thm:bound}, for $m,n$ big enough,
$$\displaystyle \Big\lceil\frac{(m+2)n}{5}\Big\rceil\leq \gamma(P_m\Box C_n).$$

\section{Wasted domination}\label{sec:wasted}

\begin{figure}[htbp]
\begin{center}
\includegraphics[width=.35\textwidth]{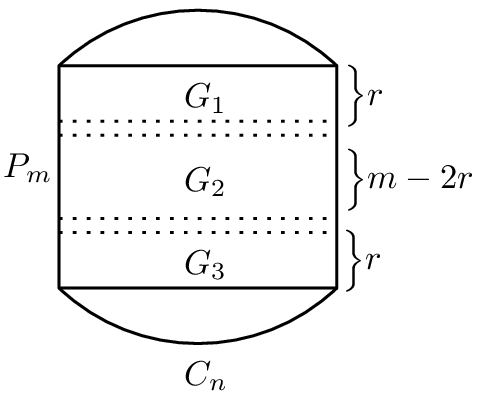}
\caption{Partition of the cylinder $P_m \Box C_n$}
\label{fig:partition}
\end{center}
\end{figure}

The wasted domination, defined in~\cite{Guichar2003} to obtain a lower bound for the domination number of grids, was later the key to computing the exact value in~\cite{Goncalves2011}.
\begin{defn}\cite{Guichar2003}
Let $G$ be a graph and let $S\subseteq V(G)$. The wasted domination of $S$ is $w(S)=5 |S| - |N[S]|$.
\end{defn}
Despite the fact that this parameter was originally defined to be used in grid graphs, it also plays an interesting role in cylinders. A vertex in a cylinder, not in the border, dominates five vertices meanwhile border vertices dominate four vertices. Therefore, for $m,n$ big enough, a dominating set of $P_m\Box C_n$ that keeps overlapping among neighborhoods of its vertices as small as possible should have a cardinal similar to $mn/5$. The ideal situation is an independent dominating set $S$ so that every vertex not in $S$ has a unique neighbor in it, that is, an efficient dominating set. Unfortunately, it is well known that the cylinder $P_m\Box C_n$ has an efficient dominating set if and only if $m=2$ and $n\equiv 0\pmod 4$ \cite{Barbosa2016}. This means that $\gamma(P_m\Box C_n)$ is greater than $mn/5$ and the wasted domination collects the information concerning the difference between both numbers. Clearly, if $S$ is a minimum dominating set of the cylinder $P_m\Box C_n$, then
$$w(S)=5 |S| - |N[S]|=5  \gamma (P_m\Box C_n) - mn,$$
or equivalently
$$\gamma (P_m\Box C_n)=\frac{w(S) + mn}{5}.$$
Suppose now that $L$ is a lower bound for $w(D)$, for every dominating set $D$ of $P_m\Box C_n$, then we obtain
\begin{equation}\label{eq:bound}
\gamma (P_m\Box C_n)\geq \frac{L + mn}{5}\cdot
\end{equation}
Our purpose is to compute a lower bound $L\leq \{ w(D)\colon \!\!D \text{\,is a dominating set of }$ $ P_m\Box C_n\}$. To this end and adapting the ideas in~\cite{Guichar2003} to the case of cylindrical graphs, we partition the cylinder $P_m\Box C_n$ into three subgraphs $G_1, G_2, G_3$ (isomorphic to cylinders with smaller paths), as is shown in Figure~\ref{fig:partition}.

Let $D$ be a dominating set of $P_m\Box C_n$. If $D_k=D\cap V(G_k)$, then
$$w(D)\geq w(D_1)+w(D_2)+w(D_3)\geq w(D_1)+w(D_3).$$

The point of the second inequality is that we expect the wasted domination to be mainly located in both borders of the cylinder, meanwhile $w(D_2)$ is zero or close to zero. Moreover, $w(D_1)$ (the same for $w(D_3)$) is computed bearing in mind that $D_1$ is a subset of $V(G_1)$ that dominates every vertex in $V(G_1)$, except perhaps vertices in the last row, and $N[D_1]$ could contain vertices not in $V(G_1)$.

Taking this construction into account, we are going to compute the minimum value of $w(R)$, where $R$ is a vertex subset of the cylinder $P_r\Box C_m$ that dominates every vertex, except perhaps vertices in the $r^{th}$ row, and $N[R]$ is computed considering that $P_r\Box C_m$ is the subgraph of $P_{r+1}\Box C_m$ consisting of the first $r$ rows. We will use the following notation.

Let $P_r\Box C_n$ be a cylinder, considered as the subgraph of the \emph{outer cylinder} $P_{r+1}\Box C_n$, consisting of the first $r$ rows.
We call \emph{inner cylinder} to $P_{r-1}\Box C_n$, the subgraph of $P_r\Box C_n$ with the first $r-1$ rows. We say that $R\subseteq V(P_r\Box C_n)$ is an \emph{almost-dominating set} of $P_r\Box C_n$ if $R$ dominates every vertex in the inner cylinder. The closed neighbor of $R$ is the  subset of the outer cylinder $N[R]=R\cup \{ u\in V(P_{r+1}\Box C_n)\colon u \text{ has a neighbor in } R\}$.

\section{The lower bound}\label{sec:bound}

The algorithm in~\cite{Guichar2003} that computes the minimum wasted domination of almost-dominating sets in grids can not be applied to cylinders, because it is specific for the Cartesian product of two paths. However, we found it inspiring, together with the ideas presented in~\cite{Klavzar1996}, in order to obtain an algorithm that plays the same role for cylinders. Having in mind that we want to improve the known lower bound (see the Inequality~\ref{eq:lower}), we will consider $P_m\Box C_n$ with $18\leq m,n$.

Let $\mathcal{P}=(\mathbb{R}\cup\{\infty\}, \min, +, \infty, 0)$ be the semi-ring of tropical numbers in the min convention (see for instance~\cite{Pin1998}), that is also called a path algebra in~\cite{Klavzar1996}. Let $\bigotimes$ be the $(\min ,+)$ matrix multiplication, that is $C=A \bigotimes B$ being the matrix where for all $i,j$, $c_{i,j}=\min \limits_{k}(a_{i,k}+b_{k,j})$.

The $(\min,+)$ product of $\alpha\in \mathbb{R}\cup\{\infty\}$ and a matrix $A$ is defined as $(\alpha\bigotimes A)_{i,j}=\alpha+a_{i,j}$. Note that $((\alpha \bigotimes A)\bigotimes B)_{i,j}=\min \limits_{k}((\alpha \bigotimes A)_{i,k}+b_{k,j})=\min \limits_{k}((\alpha+a_{i,k})+b_{k,j})=
\min \limits_{k}(\alpha+(a_{i,k}+b_{k,j}))=\alpha+\min \limits_{k}(a_{i,k}+b_{k,j})=\alpha + (A\bigotimes B)_{i,j}$. Por lo tanto $(\alpha \bigotimes A)\bigotimes B=\alpha \bigotimes(A\bigotimes B)$.

Let $G$ be a digraph with vertex set $V(G)= \{v_1, v_2, \dots, v_s\}$ together with a labeling function $\ell$ which assigns to every arc of $G$ an element of $\mathcal{P}$. A path of length $k$ in $G$ is a sequence of $k$ consecutive arcs $Q=(v_{i_0}v_{i_1})(v_{i_1}v_{i_2})\dots (v_{i_{k-1}}v_{i_k})$ and we say that $Q$ is \emph{a closed path} if $v_{i_0}=v_{i_{k}}$. The labeling can be trivially extended to paths: $\ell(Q) = \ell(v_{i_0}v_{i_1})+\ell(v_{i_1}v_{i_2})+\dots +\ell(v_{i_{k-1}}v_{i_{k}}).$

The next theorem comes from~\cite{Carre1979} and it can also be found in~\cite{Klavzar1996}. In both cases, the authors present a more general version, but we only need the following particular case related to the tropical semi-ring.

\begin{thm}\label{thm:carre}
Let $S_{ij}^k$ be the set of all paths of length $k$ from $v_i$ to $v_j$ in $G$ and let $A(G)$ be the
matrix defined by
$$A(G)_{ij} =
\left\{
\begin{array}{ll}
\ell(v_i,v_j) & \text{if }(v_i,v_j) \text{ is an arc of } G,\\
\infty & \text{otherwise.}
\end{array}
\right.
$$
If $A^k$ is the $k$-th $(\min, +)$ power of $A$, then $(A(G)^k)_{ij}=\min \{\ell(Q)\colon Q\in S_{ij}^k\}$.
\end{thm}

We will apply this result to an appropriate digraph and to this end, we need the following construction, which is a modification of the  construction for dominating sets in grid graphs described in~\cite{Spalding1998}.

Let $R\subseteq V(P_r\Box C_n)$ be an almost-dominating set and let $u\in V(P_r\Box C_n)$. We identify vertex $u$ with an element of the set $\{0,1,2\}$ according to the following rules:
\begin{itemize}
\item[ ] $u=0$ if $u\in R$;
\item[ ] $u=1$ if $u$ has at least one neighbor in $R$ in its column or in the previous one;
\item[ ] $u=2$ if $u$ has no neighbors in $R$ in its column or in the previous one.
\end{itemize}

Each column of $P_r\Box C_n$ can be seen as a word of length $r$ in the alphabet $\{0,1,2\}$, not containing the sequences $02$ nor $20$, by definition of the labeling. We denote such words as \emph{correct words} and $\mathcal{C}_r=\{p_1,p_2,\dots , p_{\alpha (r)}\}$ is the set of all correct words of length $r$. It is known that (see \cite{Goncalves2011})
$$\alpha (r)=\frac{(1+\sqrt{2})^{r+1}+(1-\sqrt{2})^{r+1}}{2}.$$

In order to identify every almost-dominating set in $P_r\Box C_n$ with an ordered list of $n$ correct words, rules to ensure that consecutive words do not violate the labeling conditions should be defined.

Let $p=(x_0,\dots , x_{r-1}), q=(y_0, \dots , y_{r-1})$ be correct words. We say that $p$ can follow $q$ if they hold the following conditions:

\begin{enumerate}
\item if $y_i=0$, then $x_i\neq 2$;
\item if $y_i=1$, then $x_i\in \{0,1,2\}$ and moreover, if $x_i=1$ then $x_{i-1}=0$ or $x_{i+1}=0$ (or both), note that in cases $i=0$ and $i=r-1$ just one between $x_{i-1}$ and $x_{i+1}$ makes sense;
\item if $y_i=2$ and $i\leq r-2$, then $x_i=0$;
\item if $y_{r-1}=2$, then $x_{r-1}\in \{0,1,2\}$ and moreover, if $x_{r-1}=1$, then $x_{r-2}=0$.
\end{enumerate}

To clarify rules 3 and 4, note that every vertex in the inner cylinder identified with $2$ has a neighbor in $R$ in the following column, in order to be dominated by $R$ (rule 3). However, vertices in the last row identified with $2$ may or may not have a neighbor in $R$, because it is not compulsory that $R$ dominates them (rule 4).

By using the labeling definition and these rules, it is clear that each almost-dominating set $R$ in $P_r\Box C_n$ has an unique associated ordered list $q_0, q_2, \dots q_{n-1}$ of $n$ correct words of length $r$ such that $q_{j+1}$ can follow $q_j$, for $j\in \{0, \dots, n-1\}$ (module $n$).

Conversely, let $q_0, q_2, \dots q_{n-1}$ be an ordered list of $n$ correct words of length $r$ such that $q_{j+1}$ can follow $q_j$, for $j\in \{0, \dots, n-1\}$ (module $n$). Denote by $q_j=(q_j^0, \dots , q_j^{r-1})$ the coordinates of $q_j$ and let $R=\{ v_{ij}\in V(P_r\Box C_m) \colon q_j^i=0 \}$. Again, by using the definition of the correct words and the rules 1,2,3 and 4, it is clear that $R$ is an almost-dominating set in $P_r\Box C_n$.

Therefore, there is an identification between almost-dominating sets in $P_r\Box C_n$ and ordered lists $q_0, q_2, \dots q_{n-1}$ of $n$ correct words of length $r$ such that $q_{j+1}$ can follow $q_j$, for $j\in \{0, \dots, n-1\}$ (module $n$). We identify the $j^{th}$ column of the cylinder with the word $q_j$, so $q_j=(v_{0,j}, \dots ,v_{r-1, j})$ and we say that $R$ \emph{starts at $q_0$}.

Let $q_j=(v_{0,j}, \dots ,v_{r-1, j})$, $q_{j+1}=(v_{0,j+1}, \dots ,v_{r-1, j+1})$ be a pair of correct words such that $q_{j+1}$ can follow $q_j$. The number $nd(q_j,q_{j+1})$ of \emph{newly dominated vertices} \cite{Guichar2003} computes how many new vertices are dominated adding $q_{j+1}$, which were not dominated before. We compute the set $\text{ND}(q_j,q_{j+1})$ of such vertices and its cardinal $nd(q_j,q_{j+1})$ as follows:

\begin{enumerate}
\item $\text{ND}=\emptyset$ and $nd=0$;
\item for each $i\in \{0,\dots, r-1\}$, if $v_{i, j+1}=0$ and $v_{i, j}=2$, then add $v_{i,j}$ to ND and add $1$ to $nd$;
\item for each $i\in \{0,\dots, r-1\}$, if $v_{i, j+1}\leq 1$ and $v_{i, j}\geq 1$, then add $v_{i, j+1}$ to ND and add $1$ to $nd$;
\item for each $i\in \{0,\dots, r-1\}$, if $v_{i, j+1}=0$, then add $v_{i, j+2}$ to ND and add $1$ to $nd$;
\item if $v_{r-1, j+1}=0$, then add $v_{r, j+1}$ to ND and add $1$ to $nd$ (note that the closed neighborhood is computed in the outer cylinder).
\end{enumerate}

\begin{lem}\label{lem:nd}
Let $R$ be an almost-dominating set of $P_r\Box C_n$ with associated ordered list $q_0, q_2, \dots q_{n-1}$. Then $|N[R]|=\sum\limits_{j=0}^{n-1}nd(q_j,q_{j+1})$.
\end{lem}

\proof
Let $R_j$ be the set of vertices of $R$ in the $j^{th}$ column, then $\text{ND}(q_j,q_{j+1})=N[R_{j+1}]\setminus (N[R_{j}]\cup N[R_{j-1}])$ and $\{\text{ND}(q_j,q_{j+1})\colon j\in \{0, \dots ,n-1\} \}$ is a partition of $N[R]$. Therefore

$$|N[R]|=\sum\limits_{j=0}^{n-1}|\text{ND}(q_j,q_{j+1})|=\sum\limits_{j=0}^{n-1}nd(q_j,q_{j+1})$$
\qed

\begin{thm}\label{th:min}
Let $G$ be the digraph with vertex set $\mathcal{C}_r=\{p_1,p_2,\dots p_{\alpha(r)}\}$ and such that there is an arc from $p_i$ to $p_j$ if $p_j$ can follow $p_i$. We denote $|p_j|$=number of zeros of $p_j$. The arc labeling is defined by
$$\ell(p_i,p_j)=5|p_j| -nd(p_i,p_j).$$
Let $A$ be the associated matrix, defined by
$$A_{ij}=
\left\{
\begin{array}{ll}
\ell(p_i,p_j) & \text{if } (p_i,p_j) \text{ is an arc, }\\
\infty & \text{otherwise.}
\end{array}
\right.
$$
Then
$$\min \limits_{1\leq i\leq \alpha(r)} (A^n)_{ii}=\min\{w(R)\colon R \text{ is an almost-dominating set of } P_r\Box C_n\}.$$
\end{thm}

\proof

A closed path $Q=(p_{i_0}p_{i_1})(p_{i_1}p_{i_2})\dots (p_{i_{n-1}}p_{i_0})$ of length $n$ is an ordered list of $n$ correct words such that $p_{i_{j+1}}$ can follow $p_{i_j}$, for $j\in \{0, \dots, n-1\}$ (module $n$). Therefore, there is a bijective correspondence between the almost-dominating sets of $P_r\Box C_n$ and the closed paths of order $n$ in $G$. Moreover, if $R$ is an almost-dominating set with associated path $Q=(p_{i_0}p_{i_2})(p_{i_2}p_{i_3})\dots (p_{i_{n-1}}p_{i_0})$, then, by Lemma~\ref{lem:nd}:

$$
\begin{array}{ll}
\ell(Q)\!& =\ell(p_{i_0}p_{i_1})+\ell(p_{i_1}p_{i_2})+\dots +\ell(p_{i_{n-1}}p_{i_{0}})\\
\!\!\!\!\!\!& = (5|p_{i_1}|\!\! - \!nd(p_{i_0},p_{i_1}))\!+\! (5 |p_{i_2}|\!\! -\!nd(p_{i_1},p_{i_2}))\!+\dots  +\!(5 |p_{i_0}|\!\! -\!nd(p_{i_{n-1}},p_{i_0}))\\
\!\!\!\!\!\!& = 5 \sum\limits_{j=0}^{n-1}|p_{i_j}| - \sum\limits_{j=0}^{n-1}nd(p_{i_j},p_{i_{j+1}})\\
\!\!\!\!\!\!& = 5|R|-|N[R]|\\
\!\!\!\!\!\!& = w(R).
\end{array}
$$

Now, using Theorem~\ref{thm:carre} we obtain that
$$
\begin{array}{ll}
(A^n)_{ii}&=\min \{\ell(Q)\colon Q \text{ is a closed path of length } n \text{ from } p_i \text{ to itself }\}\\
 & =\min \{ w(R)\!\colon\! R \text{ is an almost-dominating set of } P_r\Box C_n \text{ starting at } p_i\}.
\end{array}
$$
Finally, $\min \limits_{1\leq i\leq \alpha(r)} (A^n)_{ii}=\min \{ w(R)\colon R \text{ is an almost-dominating set of } P_r\Box C_n \}$.
\qed
\par\bigskip

We have computed the matrix $A$ and its $(\min,+)$ powers for $P_r\Box C_n$, with different values of $r$ and we have found that the case $P_{10}\Box C_n$ fits our requirements. The computation of the matrix powers has been carried out by a modification of the routine MatrixMult, available in the NVIDIA CUDA TOOLKIT \cite{nvidia}, to adapt it to the $(\min ,+)$ multiplication.

%The source code of the algorithms, in programming language C, to compute the matrix and its powers and instructions to
%generate the executable files and links to the additional libraries necessary for its compilation are available online at %https://github.com/hpcjmart/??.

We have run them in a GPU NVIDIA GeForce GTX 680. The matrix for the case $r=10$ has size $8119\times 8119$ and its memory size is about 250MB. The running time for computing each power of the matrix is $16$ seconds (it includes reading and writing the matrices in the file system).

In addition to Theorem~\ref{th:min}, we need to apply the standard recurrence argument for the $(\min,+)$ matrix multiplication (see~\cite{Goncalves2011,Guichar2003,Klavzar1996,Pavlic2013,Spalding1998}). In the case $r=10$, the computation of matrix $A$ and its powers gives that
$$(A^{120})_{i,j}=1+(A^{119})_{i,j}, \text{ for every } i,j,$$
that is, $A^{119+1}=1\bigotimes A^{119}$. We use this equality as the base case to proof by induction that $A^{n+1}=1\bigotimes A^{n}, n\geq 119$.

Assume that it is true for some $n\geq 119$, then $A^{(n+1)+1}=A^{n+1}\bigotimes A=(1\bigotimes A^{n})\bigotimes A=1\bigotimes (A^{n}\bigotimes A)=1\bigotimes A^{n+1}$, as desired. Therefore
$$(A^{n+1})_{i,j}=1+(A^{n})_{i,j}, \text{ for every } i,j \text{ and every }n\geq 119.$$
In particular, $\min \limits_{1\leq i\leq \alpha(10)} (A^{n+1})_{ii}=1+\min \limits_{1\leq i\leq \alpha(10)} (A^{n})_{ii},\text{ for every } n\geq 119.$ We have computed $A^{119}$ and we have obtained that $\min \limits_{1\leq i\leq \alpha(10)} (A^{119})_{ii}=119$.
Therefore, for $n\geq 119$:
$$
L(n)=\min \limits_{1\leq i\leq \alpha(10)} (A^n)_{ii}=n.
$$
Regarding values of $n<119$, we have individually computed each of them by using the Theorem~\ref{th:min} and the powers of the matrix $A$. We have obtained that, for $n\geq 30$:
$$L(n)=
\left\{
\begin{array}{ll}
n+1 & \text{if } n=32,33,37,38,42,43,47,48,53,58,63,\\
n & \text{otherwise.}

\end{array}
\right.
$$
The values of $L(n)$ for $n<30$ do not follow a regular formula. Finally, we can prove the announced lower bound and its immediate corollary.

\begin{thm}\label{thm:bound}
Let $P_m\Box C_n$ be a cylinder with $m\geq 20$ and $n\geq 30$, then

$$\gamma(P_m\Box C_n)\geq
\left\{
\begin{array}{ll}
\Big\lceil \displaystyle\frac{n(m+2)+2}{5}\Big\rceil& \text{if } n=32,33,37,38,42,43,47,48,53,58,63,\\
& \\
\displaystyle\Big\lceil \frac{n(m+2)}{5}\Big\rceil & \text{otherwise.}
\end{array}
\right.
$$
\end{thm}

\proof
We partition the cylinder $P_m\Box C_n$ as in Figure~\ref{fig:partition}, with $G_1$ and $G_3$ isomorphic to $P_{10}\Box C_n$. Let $D$ be a dominating set of $P_m\Box C_n$, then
$$w(D)\geq w(D_1)+w(D_3)\geq 2L(n)=L.$$
We now use this value of $L$ in the inequality $\gamma (P_m\Box C_n)\geq \frac{L + mn}{5}$, provided in Section~\ref{sec:wasted}, and we obtain
$$\gamma (P_m\Box C_n)\geq \frac{L + mn}{5}\geq
\left\{
\begin{array}{ll}
\displaystyle\frac{n(m+2)+2}{5}&\text{if } n=32,33,37,38,42,43,47,48,\\
& {\color{white} \text{if } n=\ } 53,58,63,\\
\displaystyle\frac{n(m+2)}{5} & \text{otherwise.}
\end{array}
\right.
$$
\qed

\begin{cor}
Let $P_m\Box C_n$ be a cylinder with $m\geq 20$ and $30\leq n\equiv 0\pmod 5$, then
$$\gamma(P_m\Box C_n)=\displaystyle\frac{n(m+2)}{5}\cdot$$
\end{cor}

\section{Conclusions}~\label{sec:conclusions}
In this paper we have proved that $\big\lceil\frac{(m+2)n}{5}\big\rceil\leq \gamma(P_m\Box C_n)$, for $m,n$ big enough. This new lower bound is the best one known up to now and it also provides the exact value $\gamma(P_m\Box C_{5r})=r(m+2)$. Therefore, we provide the formula for the domination number of a family of cylinders with non-bounded path length and non-bounded cycle length.

We recall here the values of $\gamma(P_m\Box C_n)$, for $16\leq m\leq 22$ and $n\geq m$, computed in~\cite{Crevals2014}, that have let us to propose this bound. Note that the particular cases in Theorem~\ref{thm:bound} with a greater lower bound, are congruent with the known values shown here. We think that there are not enough computed cases to conjecture a general formula for $\gamma(P_m\Box C_n)$, and additional work is needed in this direction.
\par\medskip

$
\gamma(P_{16}\Box C_n)=\left\{
\def\arraystretch{1.45}
\begin{array}{ll}
\big\lceil \frac{18n}{5}\big\rceil & \text{if }n\equiv 0,1,2,4\pmod 5,\\
\big\lceil \frac{18n}{5}\big\rceil +1 & \text{if }n\equiv 3\pmod 5.
\end{array}
\right.
$
\par\medskip
$
\gamma(P_{17}\Box C_n)=\left\{
\def\arraystretch{1.45}
\begin{array}{ll}
\big\lceil \frac{19n}{5}\big\rceil & \text{if }n\equiv 0,2,4\pmod 5,\\
\big\lceil \frac{19n}{5}\big\rceil +1 & \text{if }n\equiv 1,3\pmod 5.
\end{array}
\right.
$
\par\medskip
$
\gamma(P_{18}\Box C_n)=\left\{
\def\arraystretch{1.45}
\begin{array}{ll}
\big\lceil \frac{20n}{5}\big\rceil & \text{if }n\equiv 0\pmod 5,\\
\big\lceil \frac{20n}{5}\big\rceil +1 & \text{if }n\equiv 1,2,3,4\pmod 5.
\end{array}
\right.
$
\par\medskip
$
\gamma(P_{19}\Box C_n)=\left\{
\def\arraystretch{1.45}
\begin{array}{ll}
\big\lceil \frac{21n}{5}\big\rceil & \text{if }n\equiv 0,1,2\pmod 5,\\
\big\lceil \frac{21n}{5}\big\rceil +1 & \text{if }n\equiv 3,4\pmod 5.
\end{array}
\right.
$
\par\medskip
$\gamma(P_{20}\Box C_n)=\left\{
\def\arraystretch{1.45}
\begin{array}{ll}
\big\lceil \frac{22n}{5}\big\rceil & \text{if }n\equiv 0\pmod 5,\\
\big\lceil \frac{22n}{5}\big\rceil +1 & \text{if }n\equiv 1,2,3,4\pmod 5.
\end{array}
\right.
$
\par\medskip
$
\gamma(P_{21}\Box C_n)=\left\{
\def\arraystretch{1.45}
\begin{array}{ll}
\big\lceil \frac{23n}{5}\big\rceil & \text{if }n\equiv 0,2\pmod 5,\\
\big\lceil \frac{23n}{5}\big\rceil +1 & \text{if }n\equiv 1,4\pmod 5,\\
\big\lceil \frac{23n}{5}\big\rceil +2 & \text{if }n\equiv 3\pmod 5.\\
\end{array}
\right.
$
\par\medskip
$
\gamma(P_{22}\Box C_n)=\left\{
\def\arraystretch{1.45}
\begin{array}{ll}
\big\lceil \frac{24n}{5}\big\rceil & \text{if }n\equiv 0\pmod 5,\\
\big\lceil \frac{24n}{5}\big\rceil +1 & \text{if }n\equiv 1,2,4\pmod 5,\\
\big\lceil \frac{24n}{5}\big\rceil +2 & \text{if }n\equiv 3\pmod 5.\\
\end{array}
\right.
$
\par\bigskip

Finally, we would like to point out that we have used a technique based on the wasted domination, that is a concept proposed in~\cite{Guichar2003} to study the domination number in grids. We have adapted these ideas to cylinders by means of Theorem~\ref{thm:carre} from~\cite{Carre1979}, that was not used to compute the wasted domination in grids. This theorem also appears in~\cite{Klavzar1996,Pavlic2013}, in the context of the representation of the cylinders as fasciagraphs and rotagraphs, to compute exact values of the domination number of cylinders with small size. In both papers, authors provide specific constructions by means of families of subsets, being our approach different, because we use a labeling of the vertices.

In fact, our approach using the labeling to obtain the lower bound wasted domination, by means of Theorem~\ref{thm:carre}, has also been applied to compute the domination number of $P_m\Box C_n$, for $m$ small enough~\cite{Garzon}. Although these values are known, our version provides a remarkably more efficient algorithm, in terms of the size of the matrices involved in the process.

\end{document}